\documentclass[12pt]{article}%
\usepackage{sw20bams}
\usepackage{amsmath}
\usepackage{amsfonts}
\usepackage{amssymb}
\usepackage{graphicx}%
\providecommand{\U}[1]{\protect\rule{.1in}{.1in}}
\begin{document}

\title{Three Random Intercepts of a Segment}
\author{Steven Finch}
\date{February 8, 2018}
\maketitle

\begin{abstract}
We construct random triangles via uniform sampling of certain families of
lines in the plane. Two examples are given. \ The word \textquotedblleft
uniform\textquotedblright\ turns out to be vague; two competing models are
examined. \ Everything we write is well-known to experts. \ Which model is
more appropriate? \ Our hope is to engage a larger audience in answering this question.

\end{abstract}

\footnotetext{Copyright \copyright \ 2018 by Steven R. Finch. All rights
reserved.}Let $\ell$ denote a planar random line with slope $\tan(\omega)$ and
$x$-intercept $\xi$, where $\omega\sim\,$Uniform$[0,\pi]$ and $\xi\sim
\,$Uniform$[-1,1]$ are independent. Let $\ell_{1}$, $\ell_{2}$, $\ell_{3}$ be
independent copies of $\ell$. The three lines determine a compact triangle
$\Delta$ almost surely. The probability density function for the maximum angle
in $\Delta$ is \cite{Grffths}
\[
f(\alpha)=\left\{
\begin{array}
[c]{lll}%
6(3\alpha-\pi)/\pi^{2} &  & \text{if }\pi/3\leq\alpha<\pi/2,\\
6(\pi-\alpha)/\pi^{2} &  & \text{if }\pi/2\leq\alpha\leq\pi,\\
0 &  & \text{otherwise}%
\end{array}
\right.
\]
and hence the probability that $\Delta$ is obtuse is $3/4$.

A variation on the preceding is to require $\omega\sim\,$Uniform$[\pi
/4,3\pi/4]$, that is, the lines $\ell_{1}$, $\ell_{2}$, $\ell_{3}$ each have
$\left\vert \text{slope}\right\vert $ exceeding $1$. The maximum angle density
here is \cite{Finch1, BjxHys}%
\[
f(\alpha)=\left\{
\begin{array}
[c]{lll}%
24(\pi-\alpha)(2\alpha-\pi)/\pi^{3} &  & \text{if }\pi/2\leq\alpha\leq\pi,\\
0 &  & \text{otherwise.}%
\end{array}
\right.
\]
The random triangle $\Delta$ is almost surely obtuse.

Gates \cite{Gates94} examined the same two problems, for \textquotedblleft
triangles generated by uniform random lines\textquotedblright, but adopted a
different probability model than the preceding. He did not elaborate on the
quoted phrase, but referred to an earlier paper \cite{Gates93}, where it is
apparent that the density for $\omega$ should be proportional to $\sin
(\omega)$. On the one hand, his model is \textit{standard} in the sense that
the measure is invariant under rigid motions \cite{KndlMrn, Santalo}. On the
other hand, it possesses a feature that vertical lines are weighted more than
horizontal lines. This curious tradeoff raises an interesting question: which
model is more appropriate when constructing random triangles?

For the unrestricted case ($0\leq\omega<\pi$), the inclination angle density
is
\[
g(\omega)=\frac{1}{2}\sin(\omega)
\]
and the maximum angle density is consequently \cite{Gates94}
\[
f(\alpha)=\left\{
\begin{array}
[c]{lll}%
\dfrac{3}{4}\left[  (3\alpha-\pi)\cos(\alpha)+2\sin(\alpha)-2\sin
(2\alpha)+\sin(3\alpha)\right]  &  & \text{if }\pi/3\leq\alpha<\pi/2,\\
\dfrac{1}{4}\left[  3(\pi-\alpha)\cos(\alpha)+3\sin(\alpha)-2\sin
(2\alpha)\right]  &  & \text{if }\pi/2\leq\alpha\leq\pi,\\
0 &  & \text{otherwise.}%
\end{array}
\right.
\]
It follows that%
\[
\mathbb{P}\left\{  \Delta\text{ is obtuse}\right\}  =2-\frac{3\pi}%
{8}=0.8219...
\]
which is larger than $3/4$. \ For the restricted case ($\pi/4\leq\omega
\leq3\pi/4$), the inclination angle density is
\[
g(\omega)=\frac{1}{\sqrt{2}}\sin(\omega)
\]
and the maximum angle density is \cite{Gates94}%
\[
f(\alpha)=\left\{
\begin{array}
[c]{lll}%
\dfrac{1}{2}\left[  \cos(\alpha)+\sin(\alpha)+\cos(2\alpha)-2\sin
(2\alpha)\right]  &  & \text{if }\pi/2\leq\alpha\leq\pi,\\
0 &  & \text{otherwise.}%
\end{array}
\right.
\]

The expressions for $f$ when $\omega$ enjoys constant weighting are simpler
than those for $f$ when $\omega$ enjoys sinusoidal weighting. \ This statement
alone does not imply that the first model is preferable to the second model;
there are other issues to consider too. To generate random triangles according
to \cite{Gates94} is only slightly more complicated than according to
\cite{Grffths, Finch1}: if $U\sim\,$Uniform$[0,1]$, then by the inverse
CDF\ method,
\[
\omega=\arccos\left(  1-2U\right)
\]
gives inclination angles for the unrestricted case and%
\[
\omega=\arccos\left(  (1-2U)/\sqrt{2}\right)
\]
gives inclination angles for the restricted case. \ 

The most compelling argument for sinusoidal weighting is its theoretical
consistency with the Poisson line process \cite{Miles64, Solmn}. \ Let us
focus on the unrestricted case. \ By Example 20 of \cite{Baddly1}, the
inclination angles $\omega_{j}$ of the lines relative to the $x$-axis are
independent and identically distributed with density $\sin(\omega)/2$ on
$[0,\pi]$. \ In words, acute angles $\approx0^{\circ}$ and obtuse angles
$\approx180^{\circ}$ are less likely than near-right angles $\approx90^{\circ
}$. \ Vertical rain wets more than slanted rain \cite{Serra}. We quote
\cite{Baddly1}:

\begin{quotation}
... although the lines of the line process have \textquotedblleft uniformly
distributed orientations\textquotedblright\ in some sense, the angles of
incidence with any fixed axis are not uniformly distributed... the probability
of `catching' a random line in a given sampling interval of the $x$-axis
depends on the orientation of the line...
\end{quotation}

\noindent and, further, \cite{Baddly2}:

\begin{quotation}
This is a classic paradox. If you consider the random lines which
\textbf{intersect a given, fixed line}, then these random lines have angles
which are non-uniformly distributed with probability density proportional to
the sine of the incidence angle. If you consider the random lines which
\textbf{intersect a given circle} then these random lines have
uniformly-distributed orientation angles. In each case the \textbf{bold text}
describes a selection or sampling operation, and sampling operations introduce bias.
\end{quotation}

\noindent We sketch a proof of this theorem in Appendix 1. \ Proofs of the
four density formulas for $f$ are not provided here; in the following section,
we choose instead to examine only a special scenario for illustration's sake.

\section{Diagonal Line}

Let us examine the restricted case ($\pi/4\leq\omega\leq3\pi/4$), initially
with constant weighting and subsequently with sinusoidal weighting. \ We
follow \cite{BjxHys} closely. \ Let $\omega_{1}=\pi/4$, $\omega_{2}$,
$\omega_{3}$ be the inclination angles of the three lines, hence the first
line is fixed as the diagonal $y=x$. Clearly $\omega_{1}<\omega_{2}$ and
$\omega_{1}<\omega_{3}$ almost surely. \ The angles $\omega_{2}$, $\omega_{3}$
are independent and identically distributed, thus $\mathbb{P}\left\{
\omega_{2}<\omega_{3}\right\}  =1/2$. \ The triangle formed by the three lines
has angles $\omega_{2}-\omega_{1}$, $\omega_{3}-\omega_{2}$, $\pi-\omega
_{3}+\omega_{1}$. \ Since $\pi/2=\pi-3\pi/4+\pi/4\leq\pi-\omega_{3}+\omega
_{1}$, the maximum angle is obviously $\alpha=5\pi/4-\omega_{3}$. \ We have%
\begin{align*}
\mathbb{P}(\alpha &  <a)=\mathbb{P}\left\{
\begin{array}
[c]{ccc}%
5\pi/4-\omega_{3}<a & \mid & \omega_{2}<\omega_{3}%
\end{array}
\right\} \\
&  =\frac{\mathbb{P}\left\{
\begin{array}
[c]{cc}%
\omega_{3}>5\pi/4-a, & \omega_{2}<\omega_{3}%
\end{array}
\right\}  }{\mathbb{P}\left\{  \omega_{2}<\omega_{3}\right\}  }\\
&  =2\,\mathbb{P}\left\{  \omega_{3}>\max\left(
\begin{array}
[c]{cc}%
\omega_{2}, & 5\pi/4-a
\end{array}
\right)  \right\} \\
&  =2\left[
{\displaystyle\int\limits_{\pi/4}^{5\pi/4-a}}
\,%
{\displaystyle\int\limits_{5\pi/4-a}^{3\pi/4}}
g(\omega_{3})g(\omega_{2})d\omega_{3}d\omega_{2}+%
{\displaystyle\int\limits_{5\pi/4-a}^{3\pi/4}}
\,%
{\displaystyle\int\limits_{\omega_{2}}^{3\pi/4}}
g(\omega_{3})g(\omega_{2})d\omega_{3}d\omega_{2}\right]  .
\end{align*}
For $g(\omega)=2/\pi$, evaluating the double integrals yields%
\[
\mathbb{P}(\alpha<a)=\frac{(2a-\pi)(3\pi-2a)}{\pi^{2}}%
\]
and, upon differentiation,
\[%
\begin{array}
[c]{ccc}%
f(\alpha)=\dfrac{8(\pi-\alpha)}{\pi^{2}}, &  & \dfrac{\pi}{2}\leq\alpha\leq
\pi.
\end{array}
\]
For $g(\omega)=\sin(\omega)/\sqrt{2}$, evaluating the double integrals yields%
\[
\mathbb{P}(\alpha<a)=\dfrac{1}{4}\left[  2-2\cos(a)-2\sin(a)-\sin(2a)\right]
\]
and, upon differentiation,%
\[%
\begin{array}
[c]{ccc}%
f(\alpha)=\dfrac{1}{2}\left[  -\cos(\alpha)+\sin(\alpha)-\cos(2\alpha)\right]
, &  & \dfrac{\pi}{2}\leq\alpha\leq\pi.
\end{array}
\]
Moments are easily calculated; the mode is $\pi/2$ for the former and%
\[
2\arctan\left[  \frac{1}{2}\left(  -3+\sqrt{17}+\sqrt{2\left(  5-\sqrt
{17}\right)  }\right)  \right]  =1.7713...
\]
for the latter. \ Identical results apply when instead the third line is fixed
as the anti-diagonal $y=-x$.

On a personal note, I\ had intended this article to be a quick follow-up to my
2011 article on random tangents to a circle \cite{Finch2} . \ Who would have
suspected that random intercepts of a segment might be so much more hazardous
than the preceding? \ Uncovering Gates' model \cite{Gates94, Gates93}
constituted a turning point in my writing. \ This humble contribution is the
uncertain outcome of several years of hesitation and delay.

The R\ package \textit{spatstat} \cite{BRT} has planar random process
simulation capabilities. \ I\ can generate Poisson lines in a sampling window
via \textit{rpoisline} and determine their inclination angles $\omega_{j}$ via
\textit{angles.psp} (with option \textit{directed=FALSE}).\ \ An elliptical
window of eccentricity $\varepsilon\approx1$ is less likely to be hit by lines
almost parallel to the major axis than by almost perpendicular lines. \ In
contrast, for a circular window ($\varepsilon=0$), all directions are equally
likely. \ Clarifying these observations more rigorously would be worthwhile
and I\ welcome thoughts on how this should be done.

\section{Appendix 1}

The ordered pair $(\xi,\omega)$ offers one representation of a line $L$,
involving the $x$-intercept $\xi$ and inclination angle $\omega$. \ Another
representation $(p,\theta)$ where $-\infty<p<\infty$ and $0\leq\theta<\pi$,
called the Hesse normal form, involves the length $|p|$ of the perpendicular
segment from $(0,0)$ to $L$ and the orientation angle $\theta$ of this
segment. \ In the definition of a Poisson line process, it is usually assumed
that $\theta\sim\,$Uniform$[0,\pi]$. \ From%
\[
x\cos(\theta)+y\sin(\theta)=p
\]
we see that%
\[%
\begin{array}
[c]{ccc}%
p=\left\{
\begin{array}
[c]{lll}%
-\xi\sin(\omega) &  & \text{if }\omega<\pi/2,\\
\xi\sin(\omega) &  & \text{if }\omega\geq\pi/2
\end{array}
\right.  , &  & \theta=\left\{
\begin{array}
[c]{lll}%
\omega+\pi/2 &  & \text{if }\omega<\pi/2,\\
\omega-\pi/2 &  & \text{if }\omega\geq\pi/2
\end{array}
\right.
\end{array}
\]
since $\cos(\omega\pm\pi/2)=\mp\sin(\omega)$. \ At first glance, it would seem
that $\omega\sim\,$Uniform$[0,\pi]$ immediately because $\theta\sim
\,$Uniform$[0,\pi]$. \ In fact, the $2\times2$ Jacobian determinant of the
transformation $(\xi,\omega)\mapsto(p,\theta)$ is $\mp\sin(\omega)$, which
implies that the density of $\omega$ is $\sin(\omega)/2$. \ Reason for the
factor of $2$:\ both $(\xi,\omega)$ and $(-\xi,\pi-\omega)$ are mapped to the
same $(p,\theta)$. \ Details of the proof in a more general setting appear in
\cite{Solmn, Wlfwtz, Morton}.

\section{Appendix 2}

We present R\ simulation output results (ten histograms in blue) graphed
against density expressions found herein (six curves in red). \ The first four
plots correspond to the first four expressions for $f$, given without proof.
\ The next two plots correspond to those associated with the diagonal line
$y=x$ scenario. \ Analysis of other scenarios involving the vertical line
$x=0$ or the horizontal line $y=0$ are left to the reader.%
\begin{figure}[ptb]%
\centering
\includegraphics[
height=6.3354in,
width=6.3354in
]%
{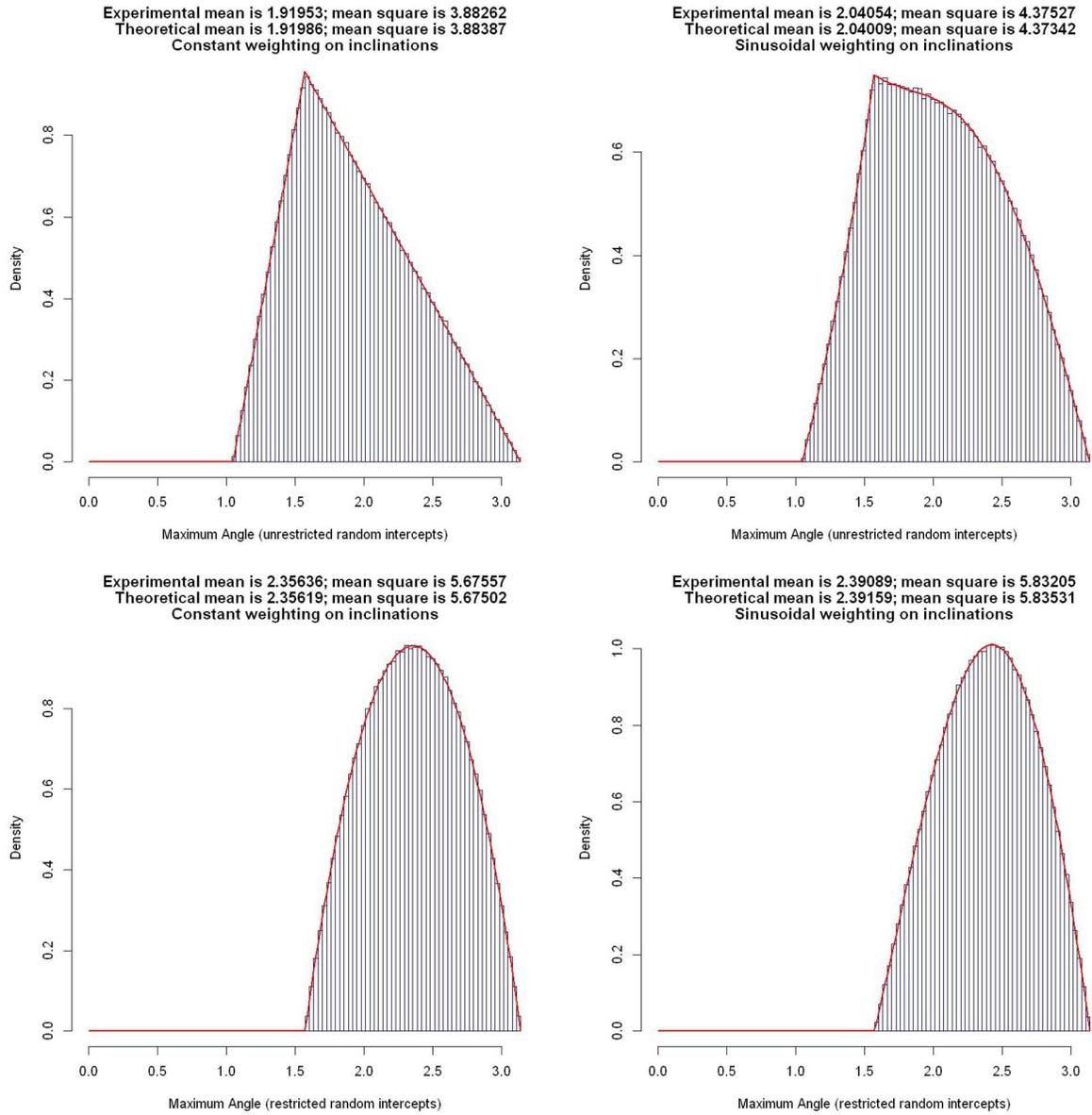}%
\caption{Top row: $f$ for unrestricted. \ Bottom row: $f$ for\ restricted.}%
\end{figure}
\begin{figure}[ptb]%
\centering
\includegraphics[
height=6.3354in,
width=4.1901in
]%
{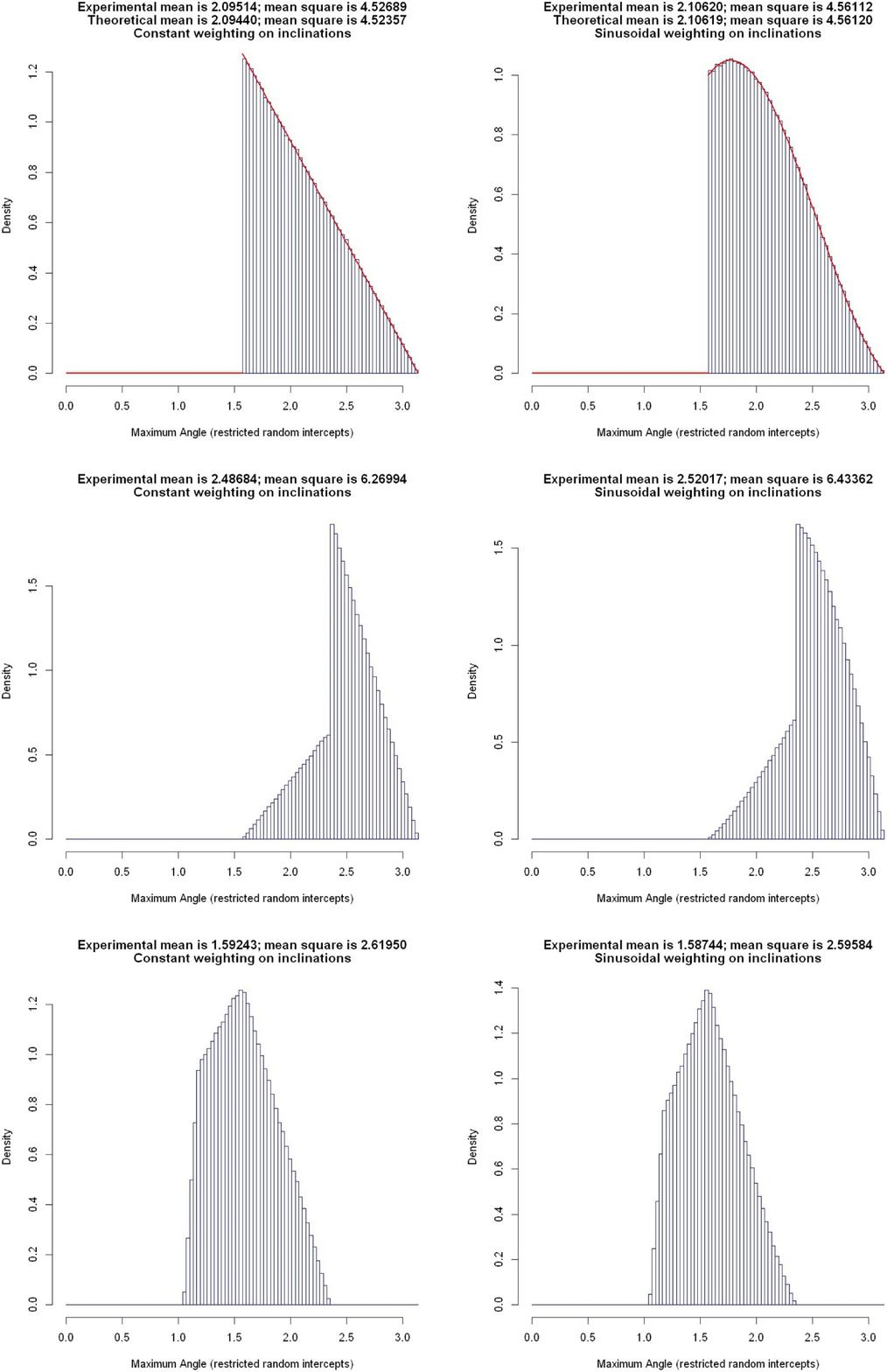}%
\caption{All restricted. \ Top row:\ diagonal $y=x$. \ Middle row: vertical
$x=0$. \ Bottom row: horizontal $y=0$.}%
\end{figure}

\section{Appendix 3}

Given a convex region $C$ in the plane, a width is the distance between a pair
of parallel $C$-supporting lines. \ Fix an inclination angle $0\leq\omega<\pi$
relative to the $x$-axis. \ A measure of all lines of angle $\omega$ hitting a
$C$-window is proportional to the corresponding width. \ For example, if $C$
is the square $[-1,1]\times\lbrack-1,1]$, we obtain a bimodal inclination
angle density \cite{Finch3}%
\[
\frac{1}{4}\max\left\{  \sqrt{1+\sin(2\omega)},\sqrt{1-\sin(2\omega)}\right\}
\]
with modes at $\pi/4$ and $3\pi/4$. \ It is easier to obtain the $\sin
(\omega)/2$ density for the interval $[-1,1]$, but harder to examine
$2\times2\sqrt{1-\varepsilon^{2}}$ rectangles of eccentricity $0<\varepsilon
<1$.

\end{document}